\documentclass[hidelinks,onefignum,onetabnum]{siamart220329} 

\usepackage{lipsum}
\usepackage{amsfonts}
\usepackage{graphicx}
\usepackage{epstopdf}
\usepackage{algorithm,algpseudocode}
\usepackage{xcolor}
\usepackage{accents}
\usepackage{multirow}
\usepackage{longtable}
\ifpdf
  \DeclareGraphicsExtensions{.eps,.pdf,.png,.jpg}
\else
  \DeclareGraphicsExtensions{.eps}
\fi

\newcommand{\bmat}[1]{\begin{bmatrix}#1\end{bmatrix}} 
\newcommand{\Lbkj}[1]{\Lbar^{(#1)}_k}
\newcommand{\Lbkjt}[1]{{\Lbar^{(#1)}_k}{}\T}

\newcommand{\Dbkj}[1]{\bar{D}^{(#1)}_k}

\newcommand{\Qkj}[1]{Q^{(#1)}_{k}}

\newcommand{\spjb}{s^+}

\newcommand{\sigp}{\sigma^+}
\newcommand{\inv}{^{-1}}
\newcommand{\subM}{_{\scriptscriptstyle M}}
\newcommand{\m}{\phantom-}
\newcommand{\words}[1]{\mgap\text{#1}\mgap}
\newcommand{\wordss}[1]{\quad\text{#1}\quad}

\newcommand{\starsymbol}{\ast}
\newcommand{\superstar}{^\starsymbol}
\newcommand{\sstar}{s\superstar}

\newcommand{\sigmastar}{\sigma\superstar}

\newcommand{\twonormd}[1]{\norm{#1}_2\drop}
\newcommand{\twonorm}[1]{\norm{#1}_2}

\newcommand{\norm}[1]{\|#1\|}
\newcommand{\Lbar}{\skew{4.3}\bar L}

\newcommand{\kp}[1]{_{k+#1}}

\newcommand{\chat}{\widehat c}
\newcommand{\ctilde}{\widetilde c}
\newcommand{\alphahat}{\skew2\widehat\alpha}
\newcommand{\etal}{et al.}  
\newcommand{\CUTEst}{{\small CUTE}st}
\newcommand{\nthiksp}{\mskip -3.5 mu}

\renewcommand{\thealgorithm}{}
\newcommand{\LDLTR}{$\mathtt{LDLtr}$}
\newcommand{\LDLtr}{$\mathtt{LDLtr}$}
\newcommand{\bfgsR}{$\mathtt{bfgsR}$}

\newcommand{\AND}{\mathop{\;\mathsf{and}\;}}
\newcommand{\OR}{\mathop{\;\mathsf{or}\;}}
\algrenewcommand{\algorithmiccomment}[1]{\hfill[#1]}
\algnewcommand{\LineComment}[1]{\State[\emph{#1}]}
\newcommand{\agap}{\hspace{1.5ex}}
\newcommand{\Astrutu}{\rule{0pt}{2.75ex}}

\newcommand{\algstrut}{\rule[-.6ex]{0pt}{2.5ex}}

\newcommand{\submax}{_{\mathrm{max}}}

\newcommand{\submin}{_{\mathrm{min}}}

\newcommand{\supplus }{^{\scriptscriptstyle +}}

\newcommand{\till}{\,{:}\,}                 

\newcommand{\mgap}{\;\;}

\newcommand{\drop}{^{\null}}
\newcommand{\kd}{_k\drop}
\newcommand{\kdp}[1]{_{k+#1}\drop}

\newcommand{\id}{_i\drop}
\newcommand{\idp}[1]{_{i+#1}\drop}

\newcommand{\T}{^{\hthinsp\mathrm{\scriptstyle T}}}

\newcommand{\invT}{^{-\mathrm{\scriptstyle T}}}

\newcommand{\hthinsp}{\mskip  1   mu}    

\newsiamremark{remark}{Remark}
\newsiamremark{hypothesis}{Hypothesis}
\crefname{hypothesis}{Hypothesis}{Hypotheses}
\newsiamthm{claim}{Claim}

\headers{$LDL\T$ Quasi-Newton}{J. J. Brust, and P. E. Gill}

\title{An $LDL\T$ Trust-Region Quasi-Newton Method \thanks{Submitted to the editors Dec. 2023.
}}

\author{Johannes J. Brust\thanks{School of Mathematical and Statistical Sciences, Arizona State University, Tempe, AZ 
  (\email{jjbrust@asu.edu}, formerly University of California San Diego).}
\and Philip E. Gill\thanks{Department of Mathematics, University of California San Diego, La Jolla, CA 
  (\email{pgill@ucsd.edu}).}
  }

\usepackage{amsopn}
\DeclareMathOperator{\diag}{diag}

\ifpdf
\hypersetup{
  pdftitle={An $LDL\T$ Quasi-Newton Trust-Region Method},
  pdfauthor={J. B. Brust, and P. E. Gill}
}
\fi

\begin{document}

\maketitle

\begin{abstract}
For quasi-Newton methods in unconstrained minimization, it is valuable to 
develop methods that are robust, i.e., methods that converge on a large number of problems.
Trust-region algorithms are often regarded to be more robust than line-search methods, however, because trust-region methods are computationally more expensive, the most popular quasi-Newton implementations use line-search methods. To fill this gap, we develop a trust-region method that updates an $LDL\T$ factorization, scales quadratically with the size of the problem, and is competitive with a conventional line-search method.
\end{abstract}



\section{Introduction}  \label{sec:intro}
Consider the unconstrained minimization problem 
\begin{equation}
\label{eq:main}
    \underset{x \in \mathbb{R}^n}{\text{ minimize}} \mgap f(x),
\end{equation}
where $f: \mathbb{R}^n \to \mathbb{R}$ is at least twice-continuously
differentiable. This problem is important for machine-learning \cite{BottouCurtisNocedal18}, model order reduction \cite{SchwerdtnerVoigt23}, financial option pricing \cite{AchdouPironneasu} and many related scientific and engineering problems. Second-order model-based methods for the problem generate an infinite sequence $\{ x_k
\}$ in which $x\kp1$ is found by minimizing a local quadratic model of $f$
based on values of the gradient $\nabla \!f(x)$ and Hessian $\nabla^2 f(x)$
at $x_k$.  Any method for unconstrained optimization must include a
globalization strategy that forces convergence from any starting
point. Broadly speaking, the two principal components of any globalization
strategy are a line search and/or the solution of a trust-region subproblem. In a conventional line-search method the local quadratic model is defined in terms of the change in variables $p_k$. Once $p_k$ has been determined, a line search is used to compute a positive scalar step length $\alpha_k$ such that $f(x_k + \alpha_k p_k)$ is sufficiently less than $f(x_k)$.  In this case the local quadratic model must be defined in terms of a positive-definite $\nabla^2 f(x_k)$ (or, a positive-definite approximation thereof) in order to ensure that the quadratic model has a bounded minimizer.  By contrast, a trust-region method is designed to give a new iterate $x_k + s_k$, where $s_k$ is a minimizer of  $\nabla \!f(x_k)\T s + \tfrac{1}{2} s\T \nabla^2 f(x_k) s$ subject to the constraint $\norm{s} \le \Delta_k$. The value of $\Delta_k$ is chosen by an iterative process designed to compute a value $f(x_k + s_k)$ that is sufficiently less than $f(x_k)$.  If the two-norm is used for the constraint $\norm{s}\le \Delta_k$ and $n$ is ``small-to-medium" in size, the standard method is due to Mor{\'e} \& Sorensen~\cite{MorS83}.  In this method the subproblem is also solved by an iterative method, with each iteration requiring the factorization of a diagonally-shifted Hessian $\nabla^2 f(x_k) + \sigma_j I$ for $\sigma_j$ a nonnegative scalar. In general, the Mor{\'e}-Sorensen method requires several factorizations to find $s_k$. However, this cost is mitigated by the fact that compared to line-search methods, trust-region methods have a stronger convergence theory and are generally more robust, i.e., they are able to solve more problems (see e.g., Dai~\cite{dai02}, Gay~\cite{Gay81}, Sorensen~\cite{Sor82}, Hebden~\cite{Heb73}, and Conn, Gould \& Toint~\cite{ConGT00a}).

If first derivatives, but not second derivatives are available, then quasi-Newton
methods can be very effective (see, e.g., Dennis \& Mor{\'e}~\cite{dennisMore77}, Gill \& Murray~\cite{gillMurray72}, Byrd, Dong \& Nocedal~\cite{byrdDongNocedal92}).  Quasi-Newton methods maintain an approximate Hessian $B_k$ (or approximate inverse Hessian $H_k$) that is modified by a low-rank update that installs the curvature information
accumulated in the step from $x_k$ to $x\kp1$.  By updating $H_k$ or a factorization of $B_k$, a quasi-Newton method can be implemented in $O(n^2)$ floating-point operations (flops) per iteration. There are infinitely many possible
modifications, but it can be argued that the most widely used quasi-Newton
method is the Broyden-Fletcher-Goldfarb-Shanno (BFGS) method (see
Broyden~\cite{Broyden70}, Fletcher~\cite{Fletcher70}, Goldfarb~\cite{Goldfarb70} and Shanno~\cite{Shanno70}), which uses a rank-two update.
This method has exhibited superior performance in a large number of
comparisons (see, Gill \& Runnoe~\cite{GilR22} for a recent survey).  In
particular, many state-of-the-art software implementations include an option
to use a quasi-Newton method, see, e.g., {\small SNOPT} \cite{snopt05}, {\small
  IPOPT} \cite{WaechterBiegler06}, Knitro \cite{knitro06} and the {\sc
  Matlab Optimization ToolBox}~\cite{MatlabOTB}.  All these implementations
use the BFGS method in conjunction with a Wolfe line-search (see, e.g.,
Mor{\'e} \& Thuente~\cite{MoreThuente94}).  In particular, if $B_k$ is positive definite and the Wolfe line-search conditions hold, then the update gives a positive-definite
matrix and the method typically exhibits a fast superlinear convergence rate. A method based on the BFGS update is the focus of Section~\ref{sec:alg}. However, we start by making no assumptions about the method used to compute $B_k$.

Although trust-region globalization methods tend to provide a more reliable
algorithm overall, the additional factorizations required at each iteration have
limited their application to quasi-Newton methods. The trust-region subproblem for the quasi-Newton case is given by
\begin{equation} \label{eq:trsub}
    \underset{\norm{s}\le \Delta_k}{\mathrm{minimize}}
           \mgap g_k\T s + \tfrac{1}{2} s\T B\kd s,
\end{equation}
where $\Delta_k$ is the positive trust-region radius, $g_k$ is the gradient $\nabla \!f(x_k)$, and $B_k$ is an $n \times n$ symmetric quasi-Newton approximation to the Hessian matrix $\nabla^2 f(x_k)$. The trust-region subproblem~(\ref{eq:trsub}) can be solved using the Mor{\'e}-Sorensen method, but the substantial cost of repeatedly factoring a shifted approximate Hessian has motivated the formulation of less expensive methods for computing $s_k$.  The most successful of these methods are based on a combination of three basic strategies.

The first strategy is to choose the norm of $s$ so that the subproblem (\ref{eq:trsub}) is easier to solve. In Gertz~\cite{gertz04} the trust-region constraint is defined in terms of the infinity-norm and the associated trust-region subproblem is solved using a quadratic programming algorithm.
An infinity-norm trust-region is also the basis of Fletcher's S$l_1$QP method for constrained optimization (see Fletcher~\cite{Fle82}).

The second strategy is to use an iterative method to solve the linear equations associated with the optimality conditions for problem (\ref{eq:trsub}). Iterative methods have the benefit of being able to compute approximate solutions of (\ref{eq:trsub}).
Steihaug~\cite{Ste83} and Toint~\cite{Toi81} apply the conjugate-gradient method to the equations $B_k s = -g_k$ but terminate the iterations if a direction of negative curvature is detected or the trust-region constraint becomes active.
 
The third strategy is to seek an approximate solution of (\ref{eq:trsub}) that lies in a  low-dimensional subspace (of dimension less than 10, say).   The dogleg method of Powell~\cite{POWELL197031} uses the subspace spanned by  the vectors $\{ g\kd, B_k^{-1}g\kd \}$. Byrd, Schnabel and Schultz~\cite{ByrdSchnabelShultz_88} propose using the subspace  based on $\{g_k, (B_k + \sigma I)^{-1}g_k  \}$ for some nonnegative $\sigma$. This extends the dogleg method to the case where $B_k$ is not positive definite. 

Other well-known iterative trust-region methods have been proposed that  do not use a quasi-Newton approximate Hessian.  These include the SSM (Hager~\cite{Hager_SSM}), GLTR (Gould, Lucidi, Roma \& Toint~\cite{GouldLucidiRomaToin_GLTR99}) and Algorithm~4 of Erway, Gill \& Griffin~\cite{ErwayGillGriffin09}.

A number of methods combine the three strategies described above. These include the methods of Brust, Marcia, Petra \& Saunders~\cite{RCR_22} and Brust, Marcia \& Petra~\cite{BMP19}.  The use of a trust-region approach in conjunction with a limited-memory approximate Hessian has been developed in Brust, Burdakov, Erway \&  Marcia~\cite{BruEM15}.

The proposed method is based on exploiting the properties of the factorization
\begin{equation} \label{eq:ldl}
  B\kd = L\kd D\kd L_k\T,
\end{equation}
where $L_k \in \mathbb{R}^{n \times n}$ is lower triangular and $D\kd \in \mathbb{R}^{n \times n}$ is diagonal. This factorization and the associated Cholesky factorization have been used extensively in the implementation of line-search quasi-Newton methods (see, e.g.,  Gill \& Murray~\cite{gillMurray72}, Fletcher \& Powell~\cite{FleP74}, Dennis \& Schnabel~\cite{DenS81}) but they are seldom used in trust-region quasi-Newton methods. In Luksan~\cite{luksan96} a factorization similar to \eqref{eq:ldl} is used for nonlinear least-squares, by introducing transformed trust-region constraints $\norm{L_k\T s} \le \Delta\kd$ in the subproblem. However, it is not applied to general minimization problems.

\subsection{Contributions}
\label{sec:contrib}
We formulate and analyze a quasi-Newton trust-region method based on exploiting the properties of the $LDL\T$ factorization.  Each iteration involves two phases. In the first phase we use a strategy similar to that proposed by Luksan~\cite{luksan96} that computes an inexpensive scalar diagonal shift for $B_k$ based on solving a trust-region subproblem with a diagonal matrix $D_k$.  In the second phase the computed shift and the factorization in \eqref{eq:ldl} are used to define an effective conjugate-gradient iteration. These steps give a quasi-Newton trust-region algorithm that is competitive with state-of-the-art line-search implementations. We note that Gould, Lucidi, Roma \& Toint~\cite{GouldLucidiRomaToin_GLTR99} have anticipated the potential of a two-phased approached, however, to the best of our knowledge, the proposed method is new.

\subsection{Notation}
\label{sec:notation}
We use Householder notation, which uses upper- and lower-case Roman letters to represent matrices and vectors, and  lower-case Greek symbols to represent scalars.
The one exception to this rule is $\Delta_k > 0$, which denotes a scalar.
The identity matrix is $I$ with dimension depending on the context. The subscript $k$ ($k\ge 0$) represents the main iteration index.  At times, an inner-iteration will be used, which is denoted by a superscript. For example, a matrix used in inner iteration $j$ of outer iteration $k$ is denoted by $Q^{(j)}_k$. The letters $D_k$ and $E_k$ are reserved for diagonal matrices and $R_k$, $L_k$ and $T_k$ denote triangular matrices.

\section{The Method} \label{sec:method}
The proposed method is based on exploiting the properties of the factorization \eqref{eq:ldl}.  First, we show how a low-rank modification of the factors can be updated in $O(n^2)$ operations.

\subsection{Updates to the factors} \label{sec:updates}
Suppose that the factorization $B\kd = L\kd D\kd L_k\T$ is available at the start of the $k$th iteration.  We make no assumptions concerning whether or not $B_k$ is positive definite but $L_k$ is assumed to be nonsingular.  We wish to compute $L_{k+1}$ and $D_{k+1}$ following a rank-one update to $B_k$. In particular, consider
\begin{equation}       \label{eq:update}
L\kd D\kd L\T_k + \alpha\kd a\kd a\T_k = L\kdp1 D\kdp1 L\kp1\T,
\end{equation}
where $a\kd \in \mathbb{R}^n$ and $\alpha_k \in \mathbb{R}$.
Let $\Lbkj{1}$ denote the $n\times(n+1)$ matrix $\bmat{L_k & a_k}$, which is lower triangular except for its last column. Similarly, let $\Dbkj{1}$ denote  the $(n+1)\times(n+1)$ diagonal matrix $\diag(D_k, \alpha_k)$. Then
\[
  L\kd D\kd L\T_k + \alpha\kd a\kd a\T_k
   = \bmat{L_k & a_k} \bmat{D_k &  \\ & \alpha_k} {\bmat{L_k & a_k}}\T
   = \Lbkj{1} \Dbkj{1} \Lbkjt{1}.
\]
A sequence of orthogonal Given's rotations $\Qkj{1} \cdots \Qkj{n}$ may be used to zero out the elements of the last column in $\Lbkj{1}$ (cf. Golub \& Van Loan~\cite{GVL4}). For $j=1$, $2$, \dots, $n$ we define
\begin{equation} \label{eq:recursion}
   \Lbkj{j+1} = \Lbkj{j}\Qkj{j},
\end{equation}
where each $\Qkj{j}$ is an $(n  + 1) \times (n+1)$ identity matrix except for four entries:
\[
\Qkj{j}(j,j) = \Qkj{j}(n+1,n+1) = \frac{\Lbkj{j}(j,j)}{\sqrt{(\Lbkj{j}(j,j))^2 + (\Lbkj{j}(j,n+1))^2}}
\]
and 
\[
  \Qkj{j}(j,n+1) = -\Qkj{j}(n+1,j)
                 =  \frac{-\Lbkj{j}(n,n+1)}{\sqrt{(\Lbkj{j}(j,j))^2 + (\Lbkj{j}(n,n+1))^2}}.
\]
The following example illustrates how $\Lbkj{1}$ is restored to triangular form:
\begin{align*}
&\bmat{ \times & & & & & \times \\ \times & \times & & & & \times \\ \times & \times & \times & & & \times \\ \times & \times & \times & \times & & \times \\ \times & \times & \times & \times & \times & \times \\} \Qkj{1} & &=
\bmat{ \times & & & & & 0\\ \times & \times & & & & \times \\ \times & \times & \times & & & \times \\ \times & \times & \times & \times & & \times \\ \times & \times & \times & \times & \times & \times \\} \\[3pt]
&\bmat{ \times & & & & & \times \\ \times & \times & & & & \times \\ \times & \times & \times & & & \times \\ \times & \times & \times & \times & & \times \\ \times & \times & \times & \times & \times & \times \\} \Qkj{1}\Qkj{2} & &=
\bmat{ \times & & & & & 0\\ \times & \times & & & & 0 \\ \times & \times & \times & & & \times \\ \times & \times & \times & \times & & \times \\ \times & \times & \times & \times & \times & \times \\} \\
& \hspace{2cm} \vdots & & \hspace{2cm} \vdots \\
&\bmat{ \times & & & & & \times \\ \times & \times & & & & \times \\ \times & \times & \times & & & \times \\ \times & \times & \times & \times & & \times \\ \times & \times & \times & \times & \times & \times \\} \Qkj{1}\Qkj{2} \cdots \Qkj{n} & &=
\bmat{ \times & & & & & 0\\ \times & \times & & & & 0 \\ \times & \times & \times & & & 0 \\ \times & \times & \times & \times & & 0 \\ \times & \times & \times & \times & \times & 0 \\}.
\end{align*}
In general, $\Lbkj{n+1} = \Lbkj{1}\Qkj{1} \cdots \Qkj{n} \equiv \bmat{L\kp1 & 0}$, and
$L\kp1$ is the first $n$ rows and columns of $\Lbkj{n+1}$. A similar recursion can be applied symmetrically to $\Dbkj{1}$ to give
\[
 \Dbkj{n+1} = (\Qkj{1}\Qkj{2} \cdots \Qkj{n})\T \Dbkj{1} \Qkj{1}\Qkj{2} \cdots \Qkj{n}                \equiv \bmat{D\kp1 & \times \\ \times & \times}.
\]
If the product  $\Qkj{1}\Qkj{2} \cdots \Qkj{n}$ is denoted by $Q_k$ then the factorization can be written as
\begin{equation}   \label{eq:ldlup}
\left\{
\begin{aligned}\mgap
    L\kd D\kd L\T_k + \alpha\kd a\kd a\T_k
      &= \Lbkj{1} \Dbkj{1} \Lbkjt{1}
       = \Lbkj{1} Q\kd Q\T_k \Dbkj{1} Q\kd Q\T_k \Lbkjt{1} \\
    &= L\kdp1 D\kdp1 L\kp1\T.
\end{aligned}
\right.
\end{equation}
If $B_k$ is not positive definite, then some of the elements of $D\kp1$ may be negative or zero and some diagonal elements of $L\kp1$ may be zero. In the latter case any offending diagonals of $L\kp1$ must be modified to give a nonsingular factor for the next iteration.

Because of the special form of each Given's rotation $\Qkj{j}$, each product in \eqref{eq:recursion} can be computed with $O(n)$ flops.  As there are $n$ total products, updating the indefinite factorization with a rank-one term requires $O(n^2)$ flops.  If more than one rank-one update is required, the  method can be applied as many times as needed. A related algorithm for updating the Cholesky factorization is given in Algorithm~C1 of Gill, Golub, Murray \& Saunders~\cite{gillGolubMurraySaunders74}.

\subsection{Computing the optimal shift}
\label{eq:shift}
Trust-region methods generate a sequence of solution estimates $\{ x_k \}$ such that 
$x_{k+1} = x_k + s_k$, where $s_k \in \mathbb{R}^n$ is a solution of the trust-region subproblem \eqref{eq:trsub}.  If the two-norm is used to define the trust region then $\sstar$ is a global minimizer of the trust-region subproblem if and only if $\twonorm{\sstar} \le \Delta_k$ and there is a scalar ``shift" $\sigmastar \ge 0$ such that
\begin{equation} \label{eq:moreSorensenOpt}
    (B_k + \sigmastar I) \sstar  = -g_k     \words{and}
      \sigmastar(\Delta_k - \twonorm{\sstar}) =  0,
\end{equation}
with $B_k + \sigmastar I$ positive semidefinite. Moreover, if $B_k + \sigmastar I$ is positive definite, then the global minimizer is unique.  Once the optimal ``shift'' $\sigmastar \ge 0$ is known, determining $\sstar$ reduces to solving the shifted linear system $(B_k + \sigmastar I) \sstar = -g_k$.  An effective algorithm due to Mor\'{e} and Sorensen \cite{MorS83} is based on using Newton's method to find a zero of the scalar-valued function $\varphi(\sigma)$ such that
\begin{equation} \label{eq:MS-iteration}
 \varphi(\sigma)
          \equiv \frac{1}{\Delta_k} - \frac{1}{\twonorm{s}},  \wordss{where}
 (B_k + \sigma I) s = -g_k.
\end{equation}
Starting with a nonnegative scalar $\sigma_0$ such that $B_k + \sigma_0 I$ is positive semidefinite, each iteration of Newton's method requires the computation of the Cholesky factorization $R\T_i R\id = B\kd + \sigma\id I$. The main computational steps of the Mor\'{e}-Sorensen method are summarized in Algorithm~1.

\begin{algorithm}[H]
\caption{$\nthiksp\bf{1}$}
\renewcommand{\thealgorithm}{}
\begin{algorithmic}[0]
\State  \Require $\sigma_0 \ge 0$; \agap
             $R_0$ such that $R_0\T R_0\drop = B\kd + \sigma_0\drop I$;\Astrutu
       \For{$i= 0\till  i\submax$}
\State    Solve $R_i\T R\id s\id = -g\kd$;
\State    Solve $R_i\T q\id = s\id$;
\State    Update $\sigma_{i+1} = \sigma_i + \dfrac{\twonorm{s_i}^2}{\twonorm{q_i}^2}                         \Big(\dfrac{\twonorm{s_i} - \Delta_k}{\Delta_k}\Big) $;
\State    Factor $R_{i+1}\T R\idp1 = B\kd + \sigma\idp1 I$;
        \EndFor;
\end{algorithmic}
\end{algorithm}


This iteration typically continues until $\twonorm{s_i} \approx \Delta_k$. Recomputing the factorization $R_{i+1}\T R\idp1$ is by far the most expensive part of the algorithm. Therefore, practical implementations typically first check whether the solution to $B_k s_0 = -g_k$ satisfies $\twonorm{s_0} \le \Delta_k$ whenever it is known that $B_k$ is positive definite to avoid this loop.

\subsection{Computing the modified shift (phase~1)} \label{sec:Modified Shift}
Since computing the optimal shift and step using Algorithm~1 is expensive, the factorization $B_k = L\kd D\kd L_k\T$ is used to compute a modified shift at a significantly reduced cost. This computation constitutes phase~1 of the proposed method.  Let $T\kd$ denote the inverse of $L_k\T$, i.e.,
\[
  L_k\T T\kd = T\kd L_k\T =  I.
\]
For any scalar $\sigma$ it holds that
\begin{equation} \label{eq:Bplussig}
  B_k + \sigma I = L\kd D\kd L_k\T + \sigma I = L\kd (D\kd + \sigma T_k\T T\kd) L_k\T.
\end{equation}
This identity can be used to modify the iteration (\ref{eq:MS-iteration}) so that expensive refactorizations are not needed. If $E_k$ is the diagonal matrix
\begin{equation}  \label{eq:Ek}
  E\kd = \diag(T_k\T T\kd),
\end{equation}
then $E_k$ can be used to approximate $T_k\T T\kd$ in the conditions (\ref{eq:moreSorensenOpt}). This gives a modified shift $\sigp$ such that
\begin{equation} \label{eq:modoptim}
  \twonorm{\spjb} \le \Delta_k, \quad  L\kd (D\kd + \sigp E\kd) L_k\T \spjb = - g\kd     \wordss{and}
      \sigp (\twonorm{\spjb} - \Delta_k) = 0,
\end{equation}
with $D\kd + \sigp E\kd$ positive semidefinite.  It is important to note that $D\kd + \sigp E\kd$ is diagonal, which allows the conditions \eqref{eq:modoptim} to be satisfied without the need for additional factorizations. The corresponding algorithm, with the initial scalar $\sigp_0 \ge 0$ is given in Algorithm~2 (details of the derivation of Algorithm~2 are given in Appendix~A).
\begin{algorithm}[H]
\caption{$\nthiksp\bf 2$\algstrut}
\begin{algorithmic}[0]
\State  \Require $\sigp_0 \ge 0$;
        \For{$i= 0\till i\submax$} \Astrutu
\State    Solve $L\kd (D\kd + \sigp_i E\kd) L_k\T \spjb_i= -g\kd$;
\State    Solve $(D\kd + \sigp_i E\kd)L_k\T q\supplus_i = -E\kd L_k\T \spjb_i$;
\State    Update $\sigp_{i+1} = \sigp_{i} - \dfrac{\twonorm{\spjb_i}^2}{ \spjb_i{}\T q\supplus_i} \Big(\dfrac{\twonorm{\spjb_i} - \Delta_k}{\Delta_k}\Big)$;
        \EndFor;
\end{algorithmic}
\end{algorithm}


Observe that Algorithm~2 requires no direct factorizations. Moreover, the solves are inexpensive because they involve only triangular or diagonal matrices. Even though the main focus of Algorithm~2 is to compute a appropriate shift $\sigp_i \ge 0$, the vector $\spjb_i$ is available as a by-product of the computation of $\sigp$.  The vector $\spjb_i$ satisfies $\twonorm{\spjb_i} = \twonorm{\sstar} = \Delta_k$ and is used to approximate $\sstar$. 

\subsection{Solving the shifted system (phase~2)}
\label{sec:phase2}
The estimate $\sigp_i$ is expected to be an overestimate to $\sigmastar$ because $\twonorm{E\kd} \le \twonorm{T_k\T T\kd} $.
Nevertheless, it contains the exact diagonal for the optimal system by \eqref{eq:Ek} and typically captures at least the right
order of magnitude. For comparison,
$\sstar$ is the solution to the shifted system $(L\kd D\kd L_k\T + \sigmastar I)\sstar = -g_k $. We use the inexpensive estimate
$\sigp_i$ to solve the related shifted system in a second phase
\begin{equation}
\label{eq:phase2}
(L\kd D\kd L_k\T + \sigp_i I) s\kd = -g\kd.
\end{equation}
The computation of an exact solution of \eqref{eq:phase2} requires a factorization of
$L\kd D\kd L_k\T + \sigp_i I$ and would be too expensive.  Instead we propose the use of an iterative solver in combination with the $LDL\T$ factors.
From \eqref{eq:Bplussig} it holds that $L\kd D\kd L_k\T + \sigp_i I = L\kd (D\kd + \sigp_i T_k\T T\kd) L_k\T$.  The conjugate-gradient (CG) method of Hestenes~\cite{HS1952} can be applied to exploit the availability of the factors:
\begin{align} \label{eq:cg1}
                         L_k h_k &= -g_k, & & (\text{triangular solve}) \\ \label{eq:cg2}
(D\kd + \sigp_i T_k\T T\kd) v\kd &= h\kd, & & (\text{conjugate-gradient solve}) \\
\label{eq:cg3}
                      L_k\T s\kd &= v\kd. & & (\text{triangular solve})
\end{align}
This vector is often close to $\sstar$ and constitutes a useful search direction.

\subsection{Backtracking the shift} \label{sec:backtrack}
As $E\kd$ is just an estimate of $T_k\T T\kd$ the computed shift $\sigp$ is usually different from $\sigmastar$. In particular, the computed $\sigp$ is often
larger than $\sigmastar$ because $\twonormd{E_k} \le \twonormd{T_k\T T\kd}$. In order to improve the accuracy of $\sigp$, a backtracking mechanism is included to allow additional trial values for $\sigp$. Specifically, the value of $\sigp$ is reduced as long as the function value decreases.  This approach is summarized in Algorithm~3.

\begin{algorithm}[H]
\caption{$\nthiksp\bf 3$\algstrut}
\begin{algorithmic}[0]
\State  \Require $\sigp_0 \ge 0$,  $\gamma_k < 1$;\Astrutu
\State    Set $i = 0$;
          \Repeat
\State      Solve $(L\kd D\kd L_k\T + \sigp_i I )s\kd = -g\kd$;
             \Comment{Use eqs. \eqref{eq:cg1}---\eqref{eq:cg3}}
\State      $\sigp_{i+1} = \gamma\kd \sigp_i$; \agap $i = i + 1$;
          \Until{$f(x_k) \le f(x_k + s_k)$};
\end{algorithmic}
\end{algorithm}

\begin{figure} \label{fig:backtrack}
\includegraphics[scale=0.55]{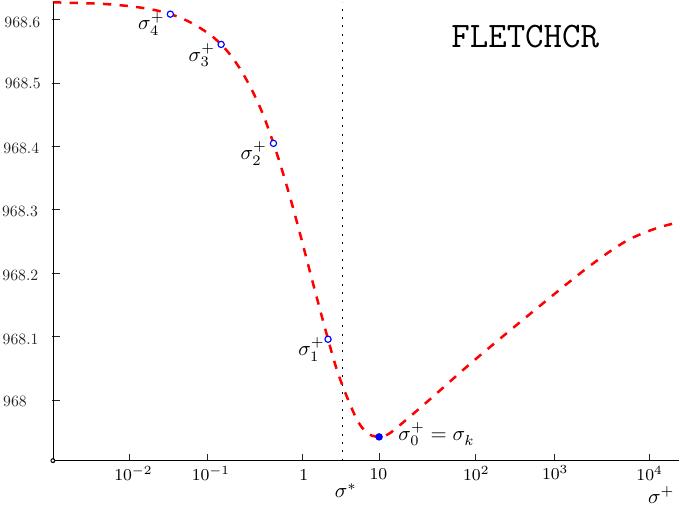}\mgap
\includegraphics[scale=0.55]{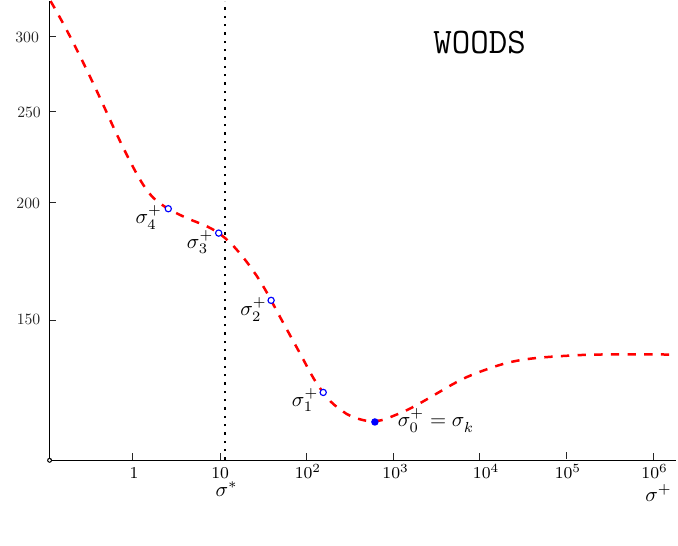}
\caption{Illustration of Alg.~3 on the two CUTEst problems \texttt{FLETCHCR} and \texttt{WOODS}. The dashed red curve represents the objective value as a function of $\sigp$. The dotted vertical line represents
the optimal shift $\sigma^*$. The blue point is the selected shift.}
\end{figure}

The search on the shift parameter is the same for both an estimated and optimal shift, i.e., the backtracking scheme could be applied if Algorithm~1 is used to solve the subproblem.  As in a backtracking line search, this strategy requires additional function evaluations. However the quality of the computed search direction is improved.  For additional efficiency, the value of $\sigp_0$ in Algorithm~3 is set to the final estimate of $\sigmastar$ computed by Algorithm~2.  Thus on initialization, $\sigp_0 \approx \sigmastar$, but as long as $f$ decreases  $\sigp_{i+1} \to 0$, and the computed steps $s_k$ become closer to the full quasi-Newton step $B_k s_k = -g_k$. Algorithm~3 
is illustrated on two problems in Fig.~\ref{fig:backtrack}, and shows that the computed $ \sigp $ can often be effective even when compared to the optimal shift.

\subsection{The quasi-Newton matrix}
\label{sec:qn}
The method can be implemented by either updating the factorization $B_k = L\kd D\kd L_k\T$ directly, or by updating its inverse. For the latter approach, recall that $T\kd$ represents the inverse of $L_k\T$, and suppose that $D\kd$ is an invertible diagonal with
inverse $G_k$. Then
\[
  D_k G_k = G_k D_k = I.
\]
and
\[
    \big(L\kd D\kd L_k\T\big){}^{-1} = L_k\invT D_k^{-1} L_k^{-1} = T\kd G\kd T_k\T.
\]
It will become evident that only $T\kd$ and $D_k$ need be stored when the inverse factorization is updated. Specifically, the proposed method generates two types of equation, with each equation associated with a particular phase. In phase~1, we solve a sequence of linear equations of the form
\begin{equation}
    L\kd \big(D\kd + \sigp_i E\kd\big) L_k\T \spjb = -g\kd.
\end{equation}
This solution can be expressed directly using only $T\kd$, $D\kd$ and $E\kd$; namely, from
\[
  h\kd  = - T_k\T g\kd, \quad w\kd = (D\kd + \sigp_i E\kd )^{-1} h\kd, \quad
  \spjb = T\kd w\kd. 
\]
Similarly, in the second phase of the method we solve systems of the form
\begin{equation*}
    L\kd (D\kd + \sigp T_k\T T\kd) L_k\T s\kd = -g\kd.
\end{equation*}
The solution of this system may also be computed using only $T\kd$ and $D\kd$. In particular,
\begin{align}
 \label{eq:phase2T1}        h\kd &= - T_k\T g\kd, & & (\text{matrix-vector multiply}) \\
 \label{eq:phase2T2}
(D\kd + \sigp_i T_k\T T\kd) v\kd &= h\kd,         & & (\text{conjugate-gradient solve}) \\
\label{eq:phase2T3}         s\kd &= T\kd v\kd. & & (\text{matrix-vector multiply})
\end{align}
 To highlight a significant difference between updating the direct factorization $L\kd D\kd L_k\T$ and updating the inverse factorization $T_k\T G\kd T\kd$ observe that the direct method computes the step in phase~2 using the equations  \eqref{eq:cg1}--\eqref{eq:cg3}.These relations depend not only on $L\kd$ and $D\kd$ but also $T\kd$. Therefore, in order to implement the direct factorization it is necessary to update $L_k$, $D\kd$ and $T\kd$. In contrast, if the inverse factorization is used, the step in phase~2 is determined by \eqref{eq:phase2T1}--\eqref{eq:phase2T3}, which depend only on $T\kd$ and $D\kd$. Therefore, updating the inverse factorization is advantageous from an implementation viewpoint because it depends only on $T\kd$ and $D\kd$. Moreover, as $D\kd$ is diagonal it is straightforward to update $T\kd$ and $G_k$, where $G_k$ is the inverse of $D_k$.  The inverse quasi-Newton matrix is denoted by $H_k$, i.e.,
\[
  B_k^{-1} = (L\kd D\kd L_k\T)^{-1}  = T\kd G\kd T_k\T \equiv H\kd.
\]
The approximate Hessian and its inverse can be positive definite or indefinite depending on the choice of updating formula. The most popular updates are defined in terms of the vectors $s_k = x\kp1 - x_k$ and $y_k = g\kp1 - g_k$. In particular the BFGS modified inverse Hessian is given by the rank-two formula
\begin{equation} \label{eq:bfgsinv}
  H\kdp1 = H\kd + \frac{y_k\T s\kd + y_k\T H\kd y\kd}{(y_k\T s\kd)^2} s\kd s_k\T - \frac1{y_k\T s\kd}(H\kd y\kd s_k\T + s\kd y_k\T H\kd),
\end{equation}
and the SR1 inverse is
\begin{equation}
\label{eq:sr1inv}
H_{k+1} = H\kd + \frac1{(y_k - H_k s_k)\T s_k}(y_k - H_k s_k) (y_k - H_k s_k)\T. 
\end{equation}
Other options are the Multipoint Symmetric Secant Matrix (MSS) update of Brust~\cite{brust2018large} and Burdakov, Mart{\'i}nez \& Pilotta~\cite{burdakov2002limited},  or the Powell-Symmetric-Broyden (PSB) update Powell~\cite{powell70} and Broyden, Dennis \& Mor{\'e}~\cite{broydenDennisMore73}.
After extensive experimentation, it was found that a quasi-Newton method based on the BFGS update \eqref{eq:bfgsinv} required the fewest function evaluations (see also Gill and Runnoe~\cite{GilR22}).  For this reason, the following discussion will focus on the
properties of a BFGS trust-region method.

Given the factorization $T\kd G\kd T_k\T$, the product $T\kdp1 G\kdp1 T\kp1\T$ is computed using a similar approach to that used in Section~\ref{sec:updates}. In particular, for the BFGS update (\ref{eq:bfgsinv}) we have
\begin{align*}  \alpha_k^{(1)} &= \frac{y_k\T s\kd + y_k\T H\kd y\kd}{(y_k\T s\kd)^2}, \quad & \alpha_k^{(2)} &=- \frac{1}{y_k\T s\kd + y_k\T H\kd y\kd},  \\
a^{(1)}_k &= s\kd + (y_k\T s\kd) \alpha^{(2)}_k H\kd y\kd, \quad & a^{(2)}_k &= H\kd y\kd.
\end{align*}
The factorization of \eqref{eq:bfgsinv} can be computed by applying \eqref{eq:ldlup} twice, i.e.,
\begin{align}
H_{k+1/2} &= T_{k+1/2}\drop G_{k+1/2}\drop T\T_{k+1/2} = T\kd G\kd T_k\T + \alpha^{(1)}_k a^{(1)}_k a^{(1)}_k{}\T \label{eq:bfgsup1}.    \\
H\kdp1   &= T\kdp1 G\kdp1 T\T_{k+1} = T_{k+1/2}\drop G_{k+1/2}\drop T\T_{k+1/2} + \alpha^{(2)}_k a^{(2)}_k a^{(2)}_k{}\T.\label{eq:bfgsup2}
\end{align}
Details of how to derive the updates are given in Appendix~B. The updates are implemented using a modification of Algorithm~C1 of Gill, Golub, Murray \& Saunders~\cite{gillGolubMurraySaunders74}.

\section{The Algorithm} \label{sec:alg}
The proposed method is given in Algorithm~4 below. 
The algorithm is a trust-region type method, with search directions being accepted when a sufficient decrease of the objective function is achieved.

\begin{algorithm}[ht]
\caption{$\nthiksp\bf 4$\algstrut}
\begin{algorithmic}[0]
\State  \Require $0 < c_1 \le c_2$, $0 < c_3 < 1 < c_4$, $0 < c_5 \le c_6 \le c_2$,
            $0 < c_7 < 1$; $0 < \gamma_0 < 1$;\Astrutu
\State  Set $T_0 = I$, $G_0 = \phi I$, $\Delta_0 > 0$,  $\epsilon >0$, $k\submax > 0$;
\State  $k=0$;
        \While{ $\epsilon \le \twonorm{g_k}  \AND k \le  k\submax$}
\State    $s\kd = - T\kd G\kd T_k\T g\kd$;
          \If{$\min_{} G_k \le 0 \OR \Delta_k < \twonorm{s_k}$}
             \If{$n\submax < n$}
\State          Compute $\spjb$, $\sigp$ from Algorithm~2; \Comment{phase~1}
             \Else
\State          Compute $\spjb$, $\sigp$ from Algorithm~1;
             \EndIf
          \EndIf
\State    Compute  $s_k$ from Algorithm~3               \Comment{phase~2}
          \If{$f(x_k + \spjb) < f(x_k + s_k)$}
\State       $s_k = \spjb$;
          \EndIf
\State    Solve $T\kd G\kd T_k\T h\kd = s\kd$ and set $\rho_k = \dfrac{f(x_k) - f(x_k + s\kd)}{s_k\T g\kd + \frac{1}{2} s_k\T h\kd}$;
          \If{$c_1 < \rho_k$}
\State       $x\kp1 = x_k + s_k$;
\State       Update $\gamma_k$;
          \Else
\State       $x\kp1 = x_k$;
          \EndIf
          \If{$c_2 < \rho_k$}
             \If{$\twonorm{s_k} \le c_3 \Delta_k$}
\State          $\Delta\kp1 = \Delta_k$;
             \Else
\State         $\Delta\kp1 = c_4\Delta_k$;   \Comment{increase trust-region radius}
             \EndIf
          \ElsIf{$ c_5 \le \rho_k \le c_6$}
\State       $\Delta\kp1 = \Delta_k$;
          \Else
\State       $\Delta\kp1 = c_7\Delta_k$;    \Comment{decrease trust-region radius}
          \EndIf
\State    Compute $G\kdp1$, $T\kdp1$ from $G\kd$, $T\kd$ using \eqref{eq:bfgsup1} and \eqref{eq:bfgsup2};                         \Comment{Update factors}
\State    $k \gets k + 1$;
        \EndWhile
\end{algorithmic}
\end{algorithm}

Two components of Algorithm~4  warrant further explanation. First, the check $n\submax < n$ branches the algorithm according to the size of the problem. As Algorithm~1 is reliable, but computationally expensive, it is used for problems that are relatively small, of the order of a hundred variables, say. For large problems a new strategy is used
to generate trial steps that estimate the shift parameter in phase~1. The trial step with the smallest objective value, becomes the next $s_k$. 
Second, if $c_1 < \rho_k$ then the step is accepted and the iterate is updated. In this case an increase or decrease of the parameter $\gamma_k$ is permitted. Specifically, subject to the limits $\gamma\submin \le \gamma_k \le  \gamma\submax$, the value of $\gamma_k$ is halved or doubled depending on the outcome of Algorithm~3. In particular, if $i=2$ then adding $\sigma\supplus_0 I$ to the quasi-Newton
matrix improved the objective, but $\gamma_k \sigma\supplus_0 I$ did not. In this case $\gamma_{k+1} = 2 \gamma_k$. On the
other hand, if $i = i\submax$ then at least every $\sigma\supplus_i$, $0 \le i \le i\submax-1$ improved the objective. In this case, $\gamma_{k+1} = \frac{1}{2}\gamma_k$. To ensure that $\gamma_{k+1}$ remains within the bounds, we set $\gamma_{k+1} \gets \max(\min(\gamma_{k+1},\gamma\submax),\gamma\submin)$. Typical values for the bounds are $\gamma\submax = \frac{1}{4} $, and $\gamma\submin = \left(\frac{1}{4}\right){}^{10}$.

\subsection{Complexity}  \label{sec:complexity}
Algorithm~4 has computational complexity of $O(n^2)$ for large $n$. To see this, note that as $T\kd$ is triangular and $G_k$ is
diagonal,  computing $s\kd = - T\kd G\kd T_k\T g\kd$ or solving $T\kd G\kd T_k\T h\kd = s\kd$ each incurs $n^2 + n \approx n^2$ multiplications. The cost of Algorithm~1 is negligible, because it is only called when
$n < n\submax$ (which is normally $n\submax = 100$). For large $n$, Algorithms~2 and 3 are used to compute
the step. Algorithm~2 is a Newton iteration for the scalar $\sigp_i$, which typically converges in 2---6 iterations. The main cost of each iteration is the solution of the triangular systems.  It is possible to achieve some savings by precomputing $L_k h_k = -g_k$ at the cost of $\frac12 n^2$ multiplications. Then, $s\supplus_i$ is obtained from $(D\kd + \sigp_i E\kd)L_k\T \spjb_i = h\kd$ in $\frac12 n^2$ multiplications. Similarly, $q\supplus_i$ is computed in $n^2$ multiplications per iteration. The overall complexity of Algorithm~2 is thus $O\left(\frac12 n^2 + i^{\mathrm{alg2}}\submax \cdot ( \frac12 n^2 + n^2 ) \right) = O\big(\frac12(3+i^{\mathrm{alg2}}\submax)n^2\big)$,
where $i^{\mathrm{alg2}}\submax$ represents the maximum iterations of Algorithm~2, which is a small integer.
Algorithm~3 implements \eqref{eq:cg1}---\eqref{eq:cg3}. The solutions of equations \eqref{eq:cg1} and \eqref{eq:cg3} are computed only once at a combined cost of $n^2$ multiplications. The conjugate gradient iteration costs $O(i^{\mathrm{cg}}\submax \cdot n^2)$
multiplications, where typically the maximum number of iterations for CG are $i^{\mathrm{cg}}\submax = 15 $.
Therefore Algorithm~3 is an $O\left( (1+ i^{\mathrm{alg3}}\submax \cdot i^{\mathrm{cg}}\submax) n^2 \right)$ computation,
where $i^{\mathrm{alg3}}\submax$ is the maximum number of iterations of Algorithm~3, typically $i^{\mathrm{alg3}}\submax =3$. 
Finally, to update the factorization in \eqref{eq:bfgsup1} and \eqref{eq:bfgsup2} the vector $H_k y_k$ is computed with
$n^2$ multiplications. It follows that the two rank-one updates of \eqref{eq:bfgsup1} and \eqref{eq:bfgsup2} are computed with $O(2n^2)$ flops. 
Combining these estimates gives the complexity of Algorithm~4 as
\begin{multline*}
O\left(\frac12\Big(4 + 3 + i^{\mathrm{alg2}}\submax + 2 + 2\cdot i^{\mathrm{alg3}}\submax \cdot i^{\mathrm{cg}}\submax + 4 \Big) n^2  \right) \\
= 
O\left(\frac12\Big(13 + i^{\mathrm{alg2}}\submax + 2\cdot i^{\mathrm{alg3}}\submax \cdot i^{\mathrm{cg}}\submax\Big) n^2  \right),
\end{multline*}
where $i^{\mathrm{alg2}}\submax$, $i^{\mathrm{alg3}}\submax$ and $i^{\mathrm{cg}}\submax$ are the maximum number of iterations for Algorithms~2, Algorithms~3, and the conjugate-gradient algorithm, respectively.  As all of these values are small constant integers, overall, Algorithm~4  is an $O(n^2)$ algorithm.
 
\subsection{Convergence}\label{sec:conv}
Algorithm~4 accepts steps that either generate a sufficient decrease or reduce the trust-region radius. From Theorem~1 of Burdakov, Gong, Yuan \& Zikrin~\cite{BurdakovLMTR16} this ensures that the trust-region algorithm converges to a stationary point of \eqref{eq:main} as long there exists a constant $\ctilde_1$ so that $\twonormd{L\kd D\kd L_k\T} \le \chat_1$, $\forall k$. This condition is equivalent to ensuring that 
$\twonorm{T\kd G\kd T_k\T} \le \chat_2$, $\forall k$ for some $\ctilde_2$.  As \eqref{eq:bfgsinv}
is positive definite when $y_k\T s\kd > 0$, $\forall k$ and $H_0 \succ 0$ we enforce these conditions  for the updates of the inverse $LDL\T$ factorization \eqref{eq:bfgsup1} and \eqref{eq:bfgsup2}. Specifically, $G_k$ and $T\kd$ are updated only if $y_k\T s\kd > 0$. Further, the initial matrix is the positive multiple of the identity $H_0\drop = T_0\drop G_0\drop T\T_0 = \phi I$ for $\phi > 0$.

\section{Implementation Details}
\label{sec:implementation}

Algorithm~4 updates the factors of $G\kd$ and $T\kd$, however the computation of $s_k$ in phase~2  using the equations \eqref{eq:phase2T1}--\eqref{eq:phase2T3} also uses the diagonal matrix $D_k$, where $D\kd = G_k\inv$. Algorithm~2 is implemented using $T\kd $ and $T_k\T$ instead of $L_k\T $ and $ L\kd$. The initial matrix is $G_0 = \phi I$, a scalar multiple of the identity where $\phi = \min\big(\max(10^{-2}, 1/\twonorm{g_0}),10^4\big)$.  The computation of  
$s_1 = x_1 - x_0$ and $y_1 = g_1 - g_0$ for the first quasi-Newton update in \eqref{eq:bfgsup1} and \eqref{eq:bfgsup2} requires the iterate $x_1$. The vector is computed using the Mor\'{e}-Thuente line-search \cite{MoreThuente94} so that $x_1 = x_0 - \alphahat_0 G_0 g_0$,
where $\alphahat_0$ is a step length that satisfies the strong Wolfe conditions. Further, we set the initial trust-region radius as $\Delta_1 = 2 \twonorm{x_1 - x_0}$. Round-off error may result in a negative diagonal element in $G_{k+1}$ when a rank-one update is made. In this case we set any negative values to their absolute values $g_j \gets | g_j |$ thereby ensuring numerical positive definiteness of $G_k$. Cancellation error can also corrupt the computation of $f(x_k) -  f(x_k + s_k)$ for determining the sufficient decrease in Algorithm~4.  Left unchecked,  the algorithm may stop making progress near a stationary point because the function values cease to provide reliable information. As a remedy, if  $f(x_k) -  f(x_k + s_k)$ is of the order of the machine precision, the sufficient decrease condition is changed to require a reduction in the gradient norm compared to $\twonorm{g_k}$. This mechanism promotes convergence to stationary points for some ill-conditioned problems.  The algorithm is implemented in Matlab and Fortran 90. 

\section{Numerical Experiments} \label{sec:experiments}
Numerical results were obtained for a large subset of the unconstrained
optimization problems from the \CUTEst{} test collection (see Bongartz
\etal{} \cite{BonCGT95} and Gould, Orban and Toint~\cite{GouOT03}).  In
particular, a problem was selected if the number of variables was of
the order of 5000 or less. The same criterion was used to set the dimension
of those problems for which the problem size can be specified. This gave a
test set of 252 problems. For comparison purposes we also give results for \bfgsR{}, which is a BFGS line-search algorithm with a line-search based on satisfying the strong Wolfe conditions. This algorithm is the state-of-the-art line-search BFGS implementation considered by Gill \& Runnoe~\cite{GilR22}. 

For assessment purposes, the cpu time,  number of iterations and number of function evaluations was recorded for each problem when
\[
   \twonorm{g_k} \le \epsilon, \words{with} \epsilon = 1 \times 10^{-4}.
\]
A limit of $k\submax = 6000$ iterations was imposed on all runs. For a given problem, if the maximum number of iterations was reached or the algorithm was unable to proceed, the data was collected if the following ``near optimal" conditions were satisfied:
\begin{equation} \label{eq:near-optimal}
 | f_k | \le |f_0| \times \epsilon\subM^{2/3}  \wordss{or} \twonorm{g_k} \le \twonorm{g_0} \times \epsilon\subM^{2/3},
\end{equation}
where $\epsilon\subM$ denotes the machine precision. Otherwise, the method was considered to have failed. Algorithm~\LDLtr{} was unable to proceed if $\Delta_k \le  10^{-22}$. Algorithm~\bfgsR{} was unable to proceed if the line search was unable to find a better point.

Details of the numerical experiments are given in the following table. An entry of ``Near opt" indicates that the method was unable to proceed but the final iterate satisfied the conditions (\ref{eq:near-optimal}).

\begin{small} 

\setlength{\tabcolsep}{2pt} 

\begin{longtable}{|ll|cccc|cccc|}
		\hline
		\multirow{2}{*}{Problem} & \multirow{2}{*}{$n$}  & 
        \multicolumn{4}{c|}{\LDLTR{} (Trust Region)} & 
		\multicolumn{4}{c|}{\bfgsR{} (Line Search)} \\
		\cline{3-10}
		&  
		& It &  Numf & Sec & Conv 
		& It &  Numf & Sec & Conv \\
		\hline
        \endfirsthead

        \multicolumn{10}{l}{(continued from the preceding page)} \\
		\hline \multirow{2}{*}{Problem} & \multirow{2}{*}{$n$}  & 
        \multicolumn{4}{c|}{\LDLTR{} (Trust Region)} & 
		\multicolumn{4}{c|}{\bfgsR{} (Line Search)} \\
		\cline{3-10}
		&  
		& It &  Numf & Sec & Conv 
		& It &  Numf & Sec & Conv \\
		\hline
        \endhead

        \hline 
        \multicolumn{10}{l}{(continued on the next page)} \\
        \endfoot
        \endlastfoot
        AKIVA&2&14&22&0.044&Opt.&15&19&0.221&Opt.\\ 
ALLINITU&4&8&10&0.045&Opt.&10&12&0.039&Opt.\\ 
ARGLINA&200&5&11&0.112&Opt.&2&4&0.017&Opt.\\ 
ARGLINB&200&{\color{orange}88}&{\color{orange}364}&{\color{orange}0.798}&{\color{orange} Near opt.}&{\color{orange}111}&{\color{orange}182}&{\color{orange}0.392}&{\color{orange} Near opt.}\\ 
ARGLINC&200&{\color{orange}91}&{\color{orange}391}&{\color{orange}0.749}&{\color{orange} Near opt.}&196&267&0.757&Opt.\\
ARGTRIGLS&200&366&2255&7.318&Opt.&206&406&0.759&Opt.\\ 
ARWHEAD&5000&7&12&5.928&Opt.&8&12&9.860&Opt.\\ 
BA-L1LS&57&60&214&0.005&Opt.&72&289&0.129&Opt.\\ 
BA-L1SPLS&57&71&260&0.004&Opt.&56&237&0.081&Opt.\\ 
BARD&3&23&25&0.022&Opt.&23&24&0.011&Opt.\\ 
BDQRTIC&5000&52&63&24.669&Opt.&34&38&40.584&Opt.\\ 
BEALE&2&17&19&0.019&Opt.&14&15&0.009&Opt.\\ 
BENNETT5LS&3&23&61&0.009&Opt.&18&26&0.016&Opt.\\ 
BIGGS6&6&44&71&0.008&Opt.&33&44&0.022&Opt.\\ 
BOX&5000&31&88&31.615&Opt.&9&13&11.476&Opt.\\ 
BOX3&3&11&13&0.014&Opt.&9&10&0.005&Opt.\\ 
BOXBODLS&2&{\color{red}---}&{\color{red}---}&{\color{red}---}&{\color{red} Max. it.}&15&37&0.018&Opt.\\ 
BOXPOWER&5000&36&64&25.848&Opt.&40&46&52.003&Opt.\\ 
BRKMCC&2&5&8&0.008&Opt.&4&7&0.004&Opt.\\ 
BROWNAL&200&5&10&0.060&Opt.&6&10&0.018&Opt.\\ 
BROWNBS&2&16&32&0.017&Opt.&{\color{orange}20}&{\color{orange}38}&{\color{orange}0.016}&{\color{orange} Near opt.}\\ 
BROWNDEN&4&19&50&0.011&Opt.&25&37&0.014&Opt.\\ 
BROYDN3DLS&5000&22&27&13.550&Opt.&23&24&28.412&Opt.\\ 
BROYDN7D&5000&360&368&206.365&Opt.&363&365&424.145&Opt.\\ 
BROYDNBDLS&5000&114&154&64.477&Opt.&53&60&62.142&Opt.\\ 
BRYBND&5000&114&154&64.472&Opt.&53&60&62.250&Opt.\\ 
CERI651ALS&7&114&197&0.044&Opt.&81&106&0.046&Opt.\\ 
CERI651BLS&7&231&530&0.004&Opt.&{\color{red}---}&{\color{red}---}&{\color{red}---}&{\color{red} Not conv.}\\ 
CERI651CLS&7&338&789&0.002&Opt.&{\color{red}---}&{\color{red}---}&{\color{red}---}&{\color{red} Not conv.}\\ 
CERI651DLS&7&{\color{red}---}&{\color{red}---}&{\color{red}---}&{\color{red} Not conv.}&206&264&0.126&Opt.\\ 
CERI651ELS&7&185&385&0.001&Opt.&100&125&0.057&Opt.\\ 
CHAINWOO&4000&537&3345&602.959&Opt.&2718&5362&1840.668&Opt.\\ 
CHNROSNB&50&167&340&0.062&Opt.&144&170&0.133&Opt.\\ 
CHNRSNBM&50&137&277&0.098&Opt.&119&147&0.075&Opt.\\ 
CHWIRUT1LS&3&28&94&0.007&Opt.&16&32&0.016&Opt.\\ 
CHWIRUT2LS&3&28&93&0.004&Opt.&15&32&0.010&Opt.\\ 
CLIFF&2&39&127&0.001&Opt.&42&60&0.022&Opt.\\ 
CLUSTERLS&2&14&16&0.031&Opt.&14&15&0.008&Opt.\\ 
COATING&134&339&834&3.140&Opt.&341&459&0.566&Opt.\\ 
COOLHANSLS&9&111&195&0.032&Opt.&172&204&0.085&Opt.\\ 
COSINE&5000&16&26&9.675&Opt.&14&22&17.041&Opt.\\ 
CRAGGLVY&5000&138&201&86.893&Opt.&279&341&327.076&Opt.\\ 
CUBE&2&51&101&0.007&Opt.&37&55&0.021&Opt.\\ 
CURLY10&5000&5118&9440&6155.342&Opt.&4258&5164&5615.681&Opt.\\ 
CURLY20&5000&4114&9947&6353.707&Opt.&3422&4536&4807.715&Opt.\\ 
CURLY30&5000&3451&9837&5549.342&Opt.&3007&4280&4136.564&Opt.\\ 
CYCLOOCFLS&4994&392&565&247.393&Opt.&378&429&504.329&Opt.\\ 
DANIWOODLS&2&12&15&0.028&Opt.&14&17&0.010&Opt.\\ 
DANWOODLS&2&115&143&0.007&Opt.&20&39&0.013&Opt.\\ 
DENSCHNA&2&10&12&0.011&Opt.&10&11&0.012&Opt.\\ 
DENSCHNB&2&7&9&0.009&Opt.&7&8&0.006&Opt.\\ 
DENSCHNC&2&14&18&0.016&Opt.&15&19&0.009&Opt.\\ 
DENSCHND&3&77&114&0.018&Opt.&65&81&0.036&Opt.\\ 
DENSCHNE&3&46&77&0.005&Opt.&30&50&0.018&Opt.\\ 
DENSCHNF&2&17&35&0.004&Opt.&10&20&0.009&Opt.\\ 
DEVGLA1&4&53&179&0.006&Opt.&33&57&0.023&Opt.\\ 
DEVGLA2&5&62&146&0.012&Opt.&40&56&0.024&Opt.\\ 
DIAMON2DLS&66&{\color{red}---}&{\color{red}---}&{\color{red}---}&{\color{red} Max. it.}&{\color{red}---}&{\color{red}---}&{\color{red}---}&{\color{red} Max. it.}\\ 
DIAMON3DLS&99&{\color{red}---}&{\color{red}---}&{\color{red}---}&{\color{red} Max. it.}&{\color{red}---}&{\color{red}---}&{\color{red}---}&{\color{red} Not conv.}\\ 
DIXMAANA&3000&22&27&5.120&Opt.&23&24&10.721&Opt.\\ 
DIXMAANB&3000&27&29&6.131&Opt.&27&29&11.912&Opt.\\ 
DIXMAANC&3000&29&31&6.515&Opt.&28&30&12.083&Opt.\\ 
DIXMAAND&3000&27&29&6.128&Opt.&28&29&11.975&Opt.\\ 
DIXMAANE&3000&1307&1312&273.515&Opt.&1308&1309&539.456&Opt.\\ 
DIXMAANF&3000&909&911&179.397&Opt.&910&911&384.064&Opt.\\ 
DIXMAANG&3000&871&873&172.185&Opt.&864&865&361.199&Opt.\\ 
DIXMAANH&3000&761&770&148.499&Opt.&735&739&303.049&Opt.\\ 
DIXMAANI&3000&{\color{red}---}&{\color{red}---}&{\color{red}---}&{\color{red} Max. it.}&{\color{red}---}&{\color{red}---}&{\color{red}---}&{\color{red} Max. it.}\\ 
DIXMAANJ&3000&1028&1030&234.822&Opt.&1029&1030&457.468&Opt.\\ 
DIXMAANK&3000&919&921&210.152&Opt.&920&921&412.102&Opt.\\ 
DIXMAANL&3000&410&412&93.968&Opt.&706&708&315.743&Opt.\\ 
DIXMAANM&3000&{\color{red}---}&{\color{red}---}&{\color{red}---}&{\color{red} Max. it.}&{\color{red}---}&{\color{red}---}&{\color{red}---}&{\color{red} Max. it.}\\ 
DIXMAANN&3000&1880&1885&370.036&Opt.&1881&1882&817.505&Opt.\\ 
DIXMAANO&3000&1508&1510&298.704&Opt.&1493&1494&664.182&Opt.\\ 
DIXMAANP&3000&1564&1566&309.012&Opt.&1566&1567&672.300&Opt.\\ 
DIXON3DQ&5000&{\color{red}---}&{\color{red}---}&{\color{red}---}&{\color{red} Max. it.}&{\color{red}---}&{\color{red}---}&{\color{red}---}&{\color{red} Max. it.}\\ 
DJTL&2&145&526&0.007&Opt.&1476&9933&2.866&Opt.\\ 
DMN15102LS&66&{\color{red}---}&{\color{red}---}&{\color{red}---}&{\color{red} Max. it.}&{\color{red}---}&{\color{red}---}&{\color{red}---}&{\color{red} Max. it.}\\ 
DMN15103LS&99&{\color{red}---}&{\color{red}---}&{\color{red}---}&{\color{red} Max. it.}&{\color{red}---}&{\color{red}---}&{\color{red}---}&{\color{red} Not conv.}\\ 
DMN15332LS&66&{\color{red}---}&{\color{red}---}&{\color{red}---}&{\color{red} Max. it.}&{\color{red}---}&{\color{red}---}&{\color{red}---}&{\color{red} Max. it.}\\ 
DMN15333LS&99&{\color{red}---}&{\color{red}---}&{\color{red}---}&{\color{red} Max. it.}&{\color{red}---}&{\color{red}---}&{\color{red}---}&{\color{red} Max. it.}\\ 
DMN37142LS&66&{\color{red}---}&{\color{red}---}&{\color{red}---}&{\color{red} Max. it.}&{\color{red}---}&{\color{red}---}&{\color{red}---}&{\color{red} Max. it.}\\ 
DMN37143LS&99&{\color{red}---}&{\color{red}---}&{\color{red}---}&{\color{red} Max. it.}&{\color{red}---}&{\color{red}---}&{\color{red}---}&{\color{red} Max. it.}\\ 
DQDRTIC&5000&28&46&13.306&Opt.&17&21&29.082&Opt.\\ 
DQRTIC&5000&799&991&537.722&Opt.&{\color{red}---}&{\color{red}---}&{\color{red}---}&{\color{red} Max. it.}\\ 
ECKERLE4LS&3&2&8&0.004&Opt.&3&5&0.003&Opt.\\ 
EDENSCH&2000&50&52&5.220&Opt.&54&56&9.811&Opt.\\ 
EG2&1000&3&6&0.182&Opt.&4&6&0.213&Opt.\\ 
EGGCRATE&2&7&12&0.005&Opt.&6&8&0.004&Opt.\\ 
EIGENALS&2550&674&1297&187.001&Opt.&665&790&201.828&Opt.\\ 
EIGENBLS&2550&5684&5717&958.874&Opt.&{\color{red}---}&{\color{red}---}&{\color{red}---}&{\color{red} Max. it.}\\ 
EIGENCLS&2652&{\color{red}---}&{\color{red}---}&{\color{red}---}&{\color{red} Max. it.}&{\color{red}---}&{\color{red}---}&{\color{red}---}&{\color{red} Max. it.}\\ 
ELATVIDU&2&13&18&0.011&Opt.&12&18&0.008&Opt.\\ 
ENGVAL1&5000&57&59&21.883&Opt.&36&37&50.878&Opt.\\ 
ENGVAL2&3&28&48&0.013&Opt.&26&31&0.016&Opt.\\ 
ENSOLS&9&23&26&0.031&Opt.&20&23&0.019&Opt.\\ 
ERRINROS&50&113&205&0.016&Opt.&247&333&0.169&Opt.\\ 
ERRINRSM&50&158&296&0.034&Opt.&329&442&0.218&Opt.\\ 
EXP2&2&10&12&0.009&Opt.&11&12&0.008&Opt.\\ 
EXPFIT&2&20&45&0.003&Opt.&10&14&0.007&Opt.\\ 
EXTROSNB&1000&125&489&7.205&Opt.&114&231&5.047&Opt.\\ 
FBRAIN3LS&6&801&2326&0.010&Opt.&1324&1735&6.021&Opt.\\ 
FLETBV3M&5000&119&150&63.915&Opt.&103&118&147.353&Opt.\\ 
FLETCBV2&5000&0&5&1.602&Opt.&0&1&0.018&Opt.\\ 
FLETCBV3&5000&8&27&8.563&Unbounded&2&11&1.769&Unbounded\\ 
FLETCHBV&5000&0&2&1.255&Unbounded&0&1&0.015&Unbounded\\ 
FLETCHCR&1000&4995&16254&368.990&Opt.&3104&5029&134.077&Opt.\\ 
FMINSRF2&4900&280&405&163.668&Opt.&243&249&317.487&Opt.\\ 
FMINSURF&4900&331&486&195.033&Opt.&289&292&376.528&Opt.\\ 
FREUROTH&5000&328&2411&887.841&Opt.&37&72&52.228&Opt.\\ 
GAUSS1LS&8&39&95&0.022&Opt.&21&32&0.014&Opt.\\ 
GAUSS2LS&8&47&98&0.003&Opt.&21&33&0.016&Opt.\\ 
GAUSS3LS&8&29&91&0.014&Opt.&24&34&0.017&Opt.\\ 
GAUSSIAN&3&1&4&0.005&Opt.&1&3&0.003&Opt.\\ 
GBRAINLS&2&8&20&0.020&Opt.&9&11&0.036&Opt.\\ 
GENHUMPS&5000&4220&14470&8968.239&Opt.&{\color{red}---}&{\color{red}---}&{\color{red}---}&{\color{red} Max. it.}\\ 
GENROSE&500&1336&4608&30.586&Opt.&838&1503&8.894&Opt.\\ 
GROWTHLS&3&1&3&0.005&Opt.&1&2&0.002&Opt.\\ 
GULF&3&54&89&0.004&Opt.&44&55&0.050&Opt.\\ 
HAHN1LS&7&275&902&0.002&Opt.&{\color{orange}107}&{\color{orange}197}&{\color{orange}0.077}&{\color{orange} Near opt.}\\ 
HAIRY&2&64&148&0.004&Opt.&18&42&0.010&Opt.\\ 
HATFLDD&3&13&15&0.012&Opt.&8&11&0.004&Opt.\\ 
HATFLDE&3&50&91&0.003&Opt.&14&15&0.009&Opt.\\ 
HATFLDFL&3&3&6&0.005&Opt.&3&5&0.004&Opt.\\ 
HATFLDFLS&3&13&24&0.005&Opt.&9&15&0.005&Opt.\\ 
HATFLDGLS&25&63&65&0.084&Opt.&66&67&0.032&Opt.\\ 
HEART6LS&6&751&1828&0.020&Opt.&3052&4151&2.706&Opt.\\ 
HEART8LS&8&2404&6450&0.001&Opt.&3133&4257&2.608&Opt.\\ 
HELIX&3&34&46&0.011&Opt.&28&35&0.011&Opt.\\ 
HIELOW&3&12&30&0.058&Opt.&13&22&0.048&Opt.\\ 
HILBERTA&2&3&8&0.003&Opt.&5&8&0.003&Opt.\\ 
HILBERTB&10&6&11&0.006&Opt.&7&8&0.005&Opt.\\ 
HIMMELBB&2&6&18&0.004&Opt.&10&19&0.006&Opt.\\ 
HIMMELBCLS&2&11&20&0.006&Opt.&7&9&0.006&Opt.\\ 
HIMMELBF&4&40&71&0.006&Opt.&36&42&0.017&Opt.\\ 
HIMMELBG&2&10&17&0.004&Opt.&5&8&0.004&Opt.\\ 
HIMMELBH&2&6&8&0.006&Opt.&5&6&0.003&Opt.\\ 
HUMPS&2&127&391&0.002&Opt.&49&124&0.030&Opt.\\ 
HYDC20LS&99&{\color{red}---}&{\color{red}---}&{\color{red}---}&{\color{red} Max. it.}&1648&1830&2.152&Opt.\\ 
INDEF&5000&7&32&6.246&Unbounded&3&12&2.871&Unbounded\\ 
INDEFM&5000&202&357&127.215&Opt.&182&213&211.380&Opt.\\ 
INTEQNELS&502&5&7&0.167&Opt.&6&7&0.092&Opt.\\ 
JENSMP&2&1&3&0.005&Opt.&1&2&0.004&Opt.\\ 
JIMACK&3549&3536&30456&9394.187&Opt.&1293&2635&1133.339&Opt.\\ 
KIRBY2LS&5&96&320&0.002&Opt.&39&60&0.021&Opt.\\ 
KOWOSB&4&26&36&0.004&Opt.&23&28&0.009&Opt.\\ 
LANCZOS1LS&6&99&166&0.003&Opt.&79&100&0.041&Opt.\\ 
LANCZOS2LS&6&99&150&0.004&Opt.&78&98&0.033&Opt.\\ 
LANCZOS3LS&6&103&163&0.006&Opt.&90&104&0.047&Opt.\\ 
LIARWHD&5000&21&31&11.977&Opt.&17&21&21.010&Opt.\\ 
LOGHAIRY&2&685&1768&0.004&Opt.&51&148&0.032&Opt.\\ 
LSC1LS&3&79&195&0.002&Opt.&55&78&0.031&Opt.\\ 
LSC2LS&3&72&92&0.001&Opt.&89&165&0.061&Opt.\\ 
LUKSAN11LS&100&1063&2380&0.127&Opt.&797&1049&1.212&Opt.\\ 
LUKSAN12LS&98&221&812&0.007&Opt.&95&196&0.110&Opt.\\ 
LUKSAN13LS&98&176&649&0.007&Opt.&44&89&0.048&Opt.\\ 
LUKSAN14LS&98&205&403&0.175&Opt.&224&273&0.236&Opt.\\ 
LUKSAN15LS&100&55&192&0.008&Opt.&68&193&0.140&Opt.\\ 
LUKSAN16LS&100&60&220&0.008&Opt.&68&249&0.106&Opt.\\ 
LUKSAN17LS&100&182&580&0.135&Opt.&421&552&0.512&Opt.\\ 
LUKSAN21LS&100&236&546&0.122&Opt.&194&260&0.318&Opt.\\ 
LUKSAN22LS&100&108&139&0.254&Opt.&115&129&0.126&Opt.\\ 
MANCINO&100&81&289&0.017&Opt.&{\color{orange}79}&{\color{orange}145}&{\color{orange}0.596}&{\color{orange} Near opt.}\\ 
MARATOSB&2&1379&3051&0.005&Opt.&1010&1421&0.678&Opt.\\ 
MEXHAT&2&40&74&0.010&Opt.&42&59&0.024&Opt.\\ 
MEYER3&3&637&1746&0.002&Opt.&{\color{orange}320}&{\color{orange}431}&{\color{orange}0.146}&{\color{orange} Near opt.}\\ 
MGH09LS&4&25&31&0.016&Opt.&16&24&0.009&Opt.\\ 
MGH10LS&3&{\color{red}---}&{\color{red}---}&{\color{red}---}&{\color{red} Not conv.}&{\color{orange}303}&{\color{orange}424}&{\color{orange}0.171}&{\color{orange} Near opt.}\\ 
MGH17LS&5&29&51&0.003&Opt.&20&32&0.010&Opt.\\ 
MISRA1ALS&2&50&130&0.007&Opt.&43&61&0.020&Opt.\\ 
MISRA1BLS&2&51&111&0.015&Opt.&34&50&0.016&Opt.\\ 
MISRA1CLS&2&26&82&0.006&Opt.&30&38&0.015&Opt.\\ 
MISRA1DLS&2&38&87&0.012&Opt.&26&35&0.015&Opt.\\ 
MNISTS0LS&494&1&3&0.317&Opt.&1&2&0.252&Opt.\\ 
MNISTS5LS&494&1&3&0.356&Opt.&1&2&0.182&Opt.\\ 
MOREBV&5000&25&165&31.345&Opt.&11&23&13.530&Opt.\\ 
MSQRTALS&1024&1860&1877&80.236&Opt.&1855&1860&81.718&Opt.\\ 
MSQRTBLS&1024&1566&1578&66.963&Opt.&1564&1569&71.230&Opt.\\ 
NCB20&5000&335&455&219.814&Opt.&208&216&259.715&Opt.\\ 
NCB20B&5000&897&908&497.100&Opt.&881&885&1189.335&Opt.\\ 
NELSONLS&3&247&774&0.001&Opt.&{\color{red}---}&{\color{red}---}&{\color{red}---}&{\color{red} Not conv.}\\ 
NONCVXU2&5000&{\color{red}---}&{\color{red}---}&{\color{red}---}&{\color{red} Max. it.}&{\color{red}---}&{\color{red}---}&{\color{red}---}&{\color{red} Max. it.}\\ 
NONCVXUN&5000&{\color{red}---}&{\color{red}---}&{\color{red}---}&{\color{red} Max. it.}&{\color{red}---}&{\color{red}---}&{\color{red}---}&{\color{red} Max. it.}\\ 
NONDIA&5000&109&770&118.168&Opt.&11&17&13.948&Opt.\\ 
NONDQUAR&5000&832&837&443.351&Opt.&831&836&999.714&Opt.\\ 
NONMSQRT&4900&{\color{red}---}&{\color{red}---}&{\color{red}---}&{\color{red} Max. it.}&{\color{red}---}&{\color{red}---}&{\color{red}---}&{\color{red} Not conv.}\\ 
OSBORNEA&5&64&101&0.011&Opt.&65&80&0.034&Opt.\\ 
OSBORNEB&11&72&134&0.032&Opt.&57&66&0.029&Opt.\\ 
OSCIGRAD&5000&486&2661&1109.817&Opt.&{\color{orange}805}&{\color{orange}3549}&{\color{orange}1189.449}&{\color{orange} Near opt.}\\ 
OSCIPATH&500&26&173&0.422&Opt.&23&47&0.292&Opt.\\ 
PALMER1C&8&53&81&0.058&Opt.&50&54&0.027&Opt.\\ 
PALMER1D&7&39&68&0.017&Opt.&37&41&0.025&Opt.\\ 
PALMER2C&8&62&86&0.008&Opt.&58&61&0.034&Opt.\\ 
PALMER3C&8&59&79&0.009&Opt.&60&63&0.028&Opt.\\ 
PALMER4C&8&60&78&0.008&Opt.&60&63&0.032&Opt.\\ 
PALMER5C&6&22&27&0.017&Opt.&24&25&0.009&Opt.\\ 
PALMER6C&8&75&94&0.007&Opt.&70&73&0.036&Opt.\\ 
PALMER7C&8&69&96&0.006&Opt.&70&73&0.037&Opt.\\ 
PALMER8C&8&70&100&0.008&Opt.&67&70&0.043&Opt.\\ 
PARKCH&15&39&163&0.221&Unbounded&51&70&3.030&Opt.\\ 
PENALTY1&1000&858&2375&46.592&Opt.&275&314&12.325&Opt.\\ 
PENALTY2&200&1371&5910&18.818&Opt.&283&971&1.161&Opt.\\ 
PENALTY3&200&{\color{orange}864}&{\color{orange}4923}&{\color{orange}40.313}&{\color{orange} Near opt.}&{\color{orange}212}&{\color{orange}522}&{\color{orange}3.435}&{\color{orange} Near opt.}\\ 
POWELLBSLS&2&191&413&0.003&Opt.&84&135&0.072&Opt.\\ 
POWELLSG&5000&44&46&26.320&Opt.&50&51&72.768&Opt.\\ 
POWER&5000&3528&8250&5353.477&Opt.&{\color{red}---}&{\color{red}---}&{\color{red}---}&{\color{red} Max. it.}\\ 
QUARTC&5000&799&991&545.104&Opt.&{\color{red}---}&{\color{red}---}&{\color{red}---}&{\color{red} Max. it.}\\ 
RAT42LS&3&1&7&0.005&Opt.&1&6&0.004&Opt.\\ 
RAT43LS&4&4&9&0.002&Opt.&5&6&0.004&Opt.\\ 
ROSENBR&2&53&101&0.001&Opt.&29&40&0.020&Opt.\\ 
ROSENBRTU&2&839&1034&0.009&Opt.&45&83&0.036&Opt.\\ 
ROSZMAN1LS&4&101&166&0.022&Opt.&24&35&0.018&Opt.\\ 
S308&2&15&17&0.011&Opt.&12&14&0.008&Opt.\\ 
SBRYBND&5000&5912&22673&15225.131&Opt.&4564&27376&8030.757&Opt.\\ 
SCHMVETT&5000&44&50&26.562&Opt.&45&47&77.080&Opt.\\ 
SCOSINE&5000&{\color{red}---}&{\color{red}---}&{\color{red}---}&{\color{red} Max. it.}&{\color{red}---}&{\color{red}---}&{\color{red}---}&{\color{red} Max. it.}\\ 
SCURLY10&5000&{\color{red}---}&{\color{red}---}&{\color{red}---}&{\color{red} Not conv.}&{\color{red}---}&{\color{red}---}&{\color{red}---}&{\color{red} Max. it.}\\ 
SCURLY20&5000&{\color{red}---}&{\color{red}---}&{\color{red}---}&{\color{red} Not conv.}&{\color{red}---}&{\color{red}---}&{\color{red}---}&{\color{red} Max. it.}\\ 
SCURLY30&5000&{\color{red}---}&{\color{red}---}&{\color{red}---}&{\color{red} Not conv.}&{\color{red}---}&{\color{red}---}&{\color{red}---}&{\color{red} Max. it.}\\ 
SENSORS&100&30&54&0.145&Opt.&24&32&0.108&Opt.\\ 
SINEVAL&2&88&209&0.003&Opt.&62&96&0.042&Opt.\\ 
SINQUAD&5000&41&77&15.391&Opt.&16&30&23.645&Opt.\\ 
SISSER&2&6&8&0.006&Opt.&4&7&0.005&Opt.\\ 
SNAIL&2&119&295&0.003&Opt.&94&133&0.049&Opt.\\ 
SPARSINE&5000&{\color{red}---}&{\color{red}---}&{\color{red}---}&{\color{red} Max. it.}&{\color{red}---}&{\color{red}---}&{\color{red}---}&{\color{red} Max. it.}\\ 
SPARSQUR&5000&230&303&149.771&Opt.&1039&1540&1365.125&Opt.\\ 
SPMSRTLS&4999&423&433&227.425&Opt.&421&426&512.698&Opt.\\ 
SROSENBR&5000&55&321&38.980&Opt.&10&15&12.296&Opt.\\ 
SSBRYBND&5000&4822&22972&12439.335&Opt.&4015&16401&4906.791&Opt.\\ 
SSCOSINE&5000&{\color{red}---}&{\color{red}---}&{\color{red}---}&{\color{red} Not conv.}&{\color{red}---}&{\color{red}---}&{\color{red}---}&{\color{red} Max. it.}\\ 
SSI&3&151&392&0.002&Opt.&1009&1418&0.649&Opt.\\ 
STRATEC&10&76&151&2.569&Opt.&47&61&1.926&Opt.\\ 
TESTQUAD&5000&{\color{red}---}&{\color{red}---}&{\color{red}---}&{\color{red} Max. it.}&954&1909&1228.821&Opt.\\ 
THURBERLS&7&110&324&0.006&Opt.&46&69&0.027&Opt.\\ 
TOINTGOR&50&154&156&0.751&Opt.&155&156&0.091&Opt.\\ 
TOINTGSS&5000&36&45&22.138&Opt.&32&34&39.242&Opt.\\ 
TOINTPSP&50&86&112&0.189&Opt.&93&111&0.069&Opt.\\ 
TOINTQOR&50&75&80&0.161&Opt.&76&77&0.043&Opt.\\ 
TQUARTIC&5000&65&299&62.010&Opt.&14&25&17.365&Opt.\\ 
TRIDIA&5000&4197&25360&10533.044&Opt.&705&1410&862.734&Opt.\\ 
VARDIM&200&47&55&0.447&Opt.&46&51&0.105&Opt.\\ 
VAREIGVL&5000&281&286&155.250&Opt.&289&291&348.909&Opt.\\ 
VESUVIALS&8&173&481&0.011&Opt.&{\color{red}---}&{\color{red}---}&{\color{red}---}&{\color{red} Max. it.}\\ 
VESUVIOLS&8&164&601&0.004&Opt.&34&70&0.051&Opt.\\ 
VESUVIOULS&8&86&294&0.018&Opt.&40&70&0.060&Opt.\\ 
VIBRBEAM&8&75&194&0.018&Opt.&75&111&0.040&Opt.\\ 
WATSON&12&61&68&0.056&Opt.&60&63&0.031&Opt.\\ 
WOODS&4000&620&3642&720.361&Opt.&531&918&361.367&Opt.\\ 
YATP1LS&4899&455&2407&875.932&Opt.&37&69&43.469&Opt.\\ 
YATP2LS&4899&86&403&63.070&Opt.&10&14&11.802&Opt.\\ 
YFITU&3&85&180&0.010&Opt.&63&82&0.029&Opt.\\ 
ZANGWIL2&2&1&3&0.007&Opt.&2&3&0.002&Opt.\\ 
\hline
\end{longtable}
\end{small}

In this experiment \LDLTR{} solved 226/252 problems to optimality,  while the line-search method solved 222/252. This indicates that the \LDLTR{} method is robust on this large subset of the unconstrained \CUTEst{} problems.

For an ``at-a-glance" comparison we provide performance profiles proposed by
Mahajan, Leyffer \& Kirches~\cite{MahajanLeyfferKirches11}, which
extend the performance profiles proposed by Dolan \&
Mor{\'e}~\cite{DolanMore02}.  In the general case with $n_p$ test problems,
performance profiles are based on values of the performance metric
\[
  \rho_s(\tau) = \frac{\text{card}\left\{ p : \pi_{p,s} \le \tau \right\}}{n_p} \words{and}
     \pi_{p,s} = \frac{t_{p,s}}{ \underset{1\le i \le S,\ i \ne s}{\text{ min}\; t_{p,i}}},
\]
where $t_{p,s}$ is the ``output'' (i.e., iterations or time) of ``solver $s$" on problem $p$, and $S$ denotes the total number of solvers for a given comparison.  When $\tau < 1$, $\rho_s(\tau)$ is an estimate of the probability that solver $s$ is faster than any other solver in $S$ by at least a factor of $1/\tau$.  For example, $\rho_s(0.25)$ is an estimate of the probability that solver $s$ is four times faster than any other solver in $S$ on a given instance.  When $\tau > 1$, $\rho_s(\tau)$ is an estimate of the probability that solver $s$ is at most $\tau$ times slower than the best-performing solver.  For example, $\rho_s(1)$ is an estimate of the probability that
solver $s$ is the fastest for a problem instance, and $\rho_s(4)$ is an estimate of the probability that solver $s$ can solve a problem at most four times slower than any other solver.

\begin{figure}[ht]
\begin{center}
\includegraphics{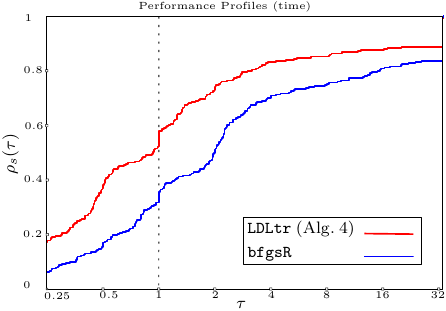}
\parbox[t]{.925\textwidth}{\caption{\label{fig:bfgsr_time}Computation-time extended performance profiles for Algorithm~4  and a strong Wolfe line-search BFGS algorithm on 252 CUTEst unconstrained problems with $n\le 5000$ variables.}}
\end{center}
\end{figure}

In Fig.~\ref{fig:bfgsr_time} we depict the performance metric
$\rho_s(\tau)$ as a function of $\tau$ for each solver $s$ (i.e., for
\bfgsR{} and \LDLTR).  A dotted vertical is used to indicate the value
$\tau = 1$.

The profiles indicate that \LDLTR{} required less overall cpu time than \bfgsR{}. The main reason for this appears to be the updating strategy of an $LDL\T$ factorization, which is implemented by updating the inverse factorization. Algorithm~4 (\LDLTR{}) is based on a modification of Algorithm~C1 of Gill, Golub, Murray \& Saunders~\cite{gillGolubMurraySaunders74}, while \bfgsR{} updates the Cholesky factors of $B_k$ using the method of Dennis \& Schnabel~\cite{DenS81}.

\begin{figure}
\begin{center}
\includegraphics{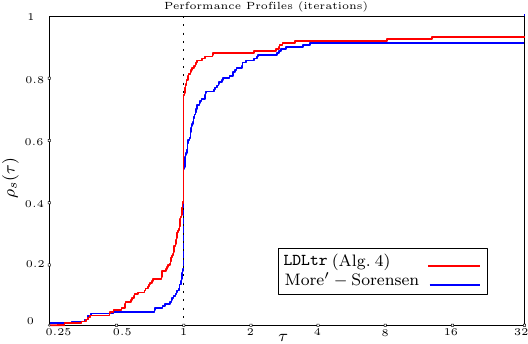}
\parbox[t]{.925\textwidth}{\caption{\label{fig:tr_comp}Extended iteration performance profiles for Algorithm~4 and the Mor\'{e}-Sorensen method on 161 CUTEst unconstrained problems with $n\le 1000$ variables.}}
\end{center}
\end{figure}

In a second experiment \LDLTR{} was compared to an implementation of the Mor\'{e}-Sorensen (MS) trust-region algorithm (\cite{MorS83}).
The MS algorithm is considered to be very robust, but requires $O(n^3)$ flops per iteration. Therefore, the problem dimension is limited to $n \le 1000$. (Note that \LDLTR{} is applicable to larger problems, but MS is not.)
The resulting test-set consisted of 161 problems. Many of the problems are relatively small, and any computational advantages in terms of time are also small. Based on the results of Fig.~\ref{fig:tr_comp}, the proposed strategy with two phases (including a branching for small and large problems) is also effective when compared to the MS algorithm.

\section{Conclusions}
\label{sec:concl}
An effective two-phase $LDL\T$ quasi-Newton trust-region algorithm has been formulated for smooth unconstrained optimization problems for which the second derivatives are not available. In the first phase, the $LDL\T$ factorization is used to compute an inexpensive estimate of the shift parameter associated with the optimality conditions for the two-norm trust-region subproblem. In the second phase, the factorization is used for a modified conjugate-gradient iteration that solves a system with the inverse approximate Hessian plus a shifted identity. Because the estimated shift parameter may be different from the optimal shift, a backtracking strategy on the shift is used to find the shift that gives the lowest function value. By updating the $LDL\T$ factorization with rank-one corrections and using two phases to generate a step, the algorithm has an overall complexity of $O(n^2)$ flops.  Numerical experiments show that the $LDL\T$ trust-region method is competitive with a strong Wolfe line-search quasi-Newton method on a subset of almost all unconstrained  problems in the \CUTEst{} test collection. The experiments indicate that the method inherits the robustness of the Mor\'{e}-Sorensen trust-region method without the computational cost.

\appendix
\section{Algorithm~2}
Algorithm~2 is based on solving for $\sigp$ and $\spjb$ in the optimality conditions \eqref{eq:modoptim}, specifically
so that $\twonorm{\spjb} = \Delta_k $. Note that, because of the first equation in \eqref{eq:modoptim}, $\spjb = \spjb(\sigp)$,  i.e., the step $\spjb$ is a function of $\sigp$. Instead of solving $\twonorm{\spjb(\sigp)} = \Delta_k$ it is better numerically 
to solve the equivalent (secular) equation
\begin{equation*}
    \phi(\sigp) = \frac{1}{\twonorm{\spjb(\sigp)}} - \frac{1}{\Delta_k} = 0.
\end{equation*}
This is a one-dimensional root finding problem in terms of $\sigp$, which can be solved with Newton's method. Starting from an initial point $\sigp_0 $, the iteration is
\begin{equation*}
    \sigp_{i+1} = \sigp_i - \frac{\phi(\sigp_i)}{\phi'(\sigp_i)}, \quad \quad i = 0,1,\ldots.
\end{equation*}
If $q\supplus = q\supplus(\sigma)$ denotes the derivative $\spjb{}' = (\spjb(\sigma))'$, then $\phi'(\sigma)$ is given by
\[
  \phi'(\sigma) = \left(\frac{1}{\twonorm{\spjb}} - \frac{1}{\Delta_k}\right)'
                = \left(\frac{1}{\spjb{}\T \spjb}\right)' 
                = - \frac{{{q\supplus}{}\T \spjb}}{\twonorm{\spjb}^3}.
\]
The quantity $q\supplus$ is computed from the equations
\begin{equation*}
    \left(L\kd (D\kd +\sigma E\kd)L_k\T \spjb = - g\kd \right){}', \quad\words{i.e.,} 
    L\kd E\kd L_k\T \spjb + L\kd(D\kd + \sigma E\kd)L_k\T q\supplus = 0,
\end{equation*}
which implies that 
\begin{equation*}
    (D\kd + \sigma E\kd) L_k\T q\supplus = - E\kd L_k\T \spjb.
\end{equation*}
The Newton correction at $\sigma$ is then
\begin{equation*}
   \frac{\phi(\sigma)}{\phi'(\sigma)}
     = \frac{\phi}{\phi'}
     = - \frac{\left( 1 / \twonorm{\spjb} - 1 / \Delta_k  \right)}{ \spjb{}\T q\supplus /          \twonorm{\spjb}^3}
     = - \frac{\twonorm{\spjb}^2}{\spjb{}\T q\supplus} \Big(\frac{\Delta_k -       \twonorm{\spjb} }{\Delta_k}\Big).
\end{equation*}
It follows that at $\sigma = \sigma_i$, we have
\begin{equation*}
    \sigma_{i+1} = \sigma_i - \frac{\twonorm{\spjb}^2}{ \spjb{}\T q\supplus} \Big(\frac{\twonorm{\spjb} - \Delta_k}{\Delta_k}\Big),
\end{equation*}
which completes the derivation of the quantities used in Algorithm~2.

\section{Quasi-Newton $LDL\T$ Updates}
In order to apply the $LDL\T$ updating strategy to \eqref{eq:bfgsinv} we reformulate the rank-2 update in \eqref{eq:bfgsinv} to conform to \eqref{eq:update}.  Specifically, let
\begin{equation*}
    \beta_1 = \frac{y_k\T s\kd + y_k\T H\kd y\kd}{(y_k\T s\kd)^2}, \qquad
    \beta_2 = \frac{1}{y_k\T s\kd},
\end{equation*}
and note that 
\begin{multline*}
    \frac{ (y_k\T s\kd + y_k\T H\kd y\kd) s\kd}{(y_k\T s\kd)^2}\,s\kd s_k\T
      - \frac1{ y_k\T s\kd}\big(H\kd y\kd s_k\T + s\kd y_k\T H\kd\big)  \\
      = \bmat{ s_k & H_k y_k}
        \bmat{ \m\beta_1 & - \beta_2 \\ -\beta_2 & \m0}
        \bmat{ s_k & H_k y_k}\T.
\end{multline*}
The $2\times 2$ matrix can be factored as
\begin{equation*}
  \bmat{ \m\beta_1 & - \beta_2 \\[5pt]
          -\beta_2 & \m  0}
   = \bmat{1              & 0 \\[5pt]
           l_{21}         & 1   } 
     \bmat{\alpha^{(1)}_k &   \\ & \alpha^{(2)}_k}
     \bmat{1 & l_{21} \\[5pt] 0 & 1}.
\end{equation*}
where
\begin{equation*}
    \alpha^{(1)}_k = \beta_1, \quad l_{21} = - \frac{\beta_2}{\beta_1}, \quad \text { and} \quad \alpha^{(2)}_k = - \frac{\beta^2_2}{\beta_1}.
\end{equation*}
If $a^{(1)}_k$ and $a^{(2)}_k$ denote the quantities
\begin{equation*}
a^{(1)}_k = \bmat{ s_k & H_k y_k} \bmat{ 1 \\ l_{21}} = s\kd + (y_k\T s\kd) \alpha^{(2)}_k H\kd y\kd, \wordss{and} 
a^{(2)}_k = H\kd y\kd,
\end{equation*}
then 
\begin{equation*}
\alpha^{(1)}_k a^{(1)}_k a^{(1)}_k{}\T + \alpha^{(2)}_k a^{(2)}_k a^{(2)}_k{}\T
 = \frac{ (y_k\T s\kd + y_k\T H\kd y\kd)}{(y_k\T s\kd)^2} s\kd s_k\T
   - \frac1{ y_k\T s\kd} \big( H\kd y\kd s_k\T + s\kd y_k\T H\kd\big).
\end{equation*}
which are the updates used in \eqref{eq:bfgsup1} and \eqref{eq:bfgsup2}.





\bibliographystyle{siamplain}
\bibliography{references}
\end{document}